\documentclass{amsart}[12pt]

\usepackage[cp866]{inputenc}\usepackage[english]{babel}
\usepackage{latexsym}\usepackage{amsfonts} \usepackage{amssymb}

\newtheorem{theorem}{Theorem}

\newtheorem{prop}{Proposition}

\newtheorem{remark}{Remark}
\newtheorem{definition}{Definition}
\newtheorem{lemma}{Lemma}

\newcommand{\Z}{{\mathbb Z}}

\newcommand{\R}{{\mathbb R}}
\newcommand{\CC}{{\mathbb C}}
\begin{document}

\title{There are no algebraically integrable ovals \\ in even-dimensional
spaces}

\author{V.A.~Vassiliev}

\date{2013}

\address{Steklov Mathematical Institute and National Research University
Higher School of Economics, Moscow, Russia}

\email{vva@mi.ras.ru}

\thanks{Supported by RFBR grant 13-1-00383 and the program ``Leading scientific
schools'', grant No. NSh-4850.2012.1}

\dedicatory{To the memory of Egbert Brieskorn}



\begin{abstract}
We prove that there are no bounded domains with smooth boundaries in
even-dimensional Euclidean spaces, such that the volumes cut off from them by
affine hyperplanes depend algebraically on these hyperplanes. For convex ovals
in ${\R}^2$, this is the Newton's Lemma XXVIII from \cite{Newton}, see also
\cite{comment}, \cite{Kvant}, \cite{Notices}.
\end{abstract}

\maketitle

\section{Introduction}

\unitlength 0.7mm \linethickness{0.6pt}

Sir Isaac Newton has proved (see \cite{Newton}, Lemma XXVIII of Book 1) that
the areas cut off by different lines from a convex bounded domain with
infinitely smooth boundary in ${\R^2}$ never define an algebraic function on
the space of lines\footnote{``There is no oval figure whose area, cut off by
right lines at pleasure, can be universally found by means of equations of any
number of finite terms and dimensions'', originally ``Nulla extat figura Ovalis
cujus area, rectis pro lubitu abscissam possit per aequationes numero
terminorum ac dimensionum finitas generaliter inveniri''}.

This fact contrasts to the Archimedes' theorem on sphere sections, implying
that the volumes cut off by planes from a ball in ${\R^3}$ depend algebraically
on these planes; it is easy to check that the last fact holds also for
arbitrary ellipsoids in odd-dimensional spaces.

In 1987, in connection with the 300th anniversary of the Newton's Book,
V.I.~Arnold has asked whether his result is true in the case of other
dimensions and general domains with smooth boundaries, see \cite{Arprob},
problems 1987-14, 1988-13 and 1990-27.

In 1988, solving this problem, I have extended the Newton's result to convex
domains in even-dimensional spaces and to arbitrary bounded domains with smooth
boundaries in ${\R}^2$, see \cite{Notices}, \cite{APLT}. In the present paper,
the same statement is proved for arbitrary bounded domains with smooth
boundaries in even-dimensional spaces. The proof is based on Picard-Lefschetz
theory and elementary facts on finite reflection groups.

\subsection{Definitions and main theorem}

Denote by $P_n$ the space of all affine hyperplanes in $\R^n$. It can be
considered as $\R P^n$ with one point removed, in particular as an algebraic
manifold. Given a compact domain $D \subset \R^n$, the corresponding two-valued
{\em volume function} $P_n \to \R$ associates with any hyperplane $L \in P_n$
the volumes of both parts cut off from $D$ by this hyperplane. The domain $D$
is called {\em algebraically integrable} if this function is algebraic, i.e.
there exists a non-trivial polynomial $\Phi$ in $n+2$ variables such that for
any real numbers $a_1, \dots, a_n, b$ and any value $V_i$, $i=1,2$, of the
volume function on the plane $L \in P_n$ defined by the equation $a_1 x_1 +
\dots + a_n x_n = b$, we have $\Phi(a_1, \dots , a_n, b, V_i)=0$. We are going
to prove the following fact.

\begin{theorem}
\label{mainth} If $n$ is even, then there is no algebraically integrable
bounded domain with $C^\infty$-smooth boundary in $\R^n$.
\end{theorem}

\begin{remark} \rm
1. For ovals (i.e. convex  bounded domains) in $\R^2$ this theorem was proved
in \cite{Newton}, see also \cite{comment}, \cite{Notices}.

2. If $D$ is a bounded algebraically integrable domain with $C^\infty$-smooth
boundary in $\R^n$, then by projective duality and Tarski--Seidenberg lemma
this boundary $\partial D$ is semialgebraic, see \cite{Kvant}, \cite{Notices}.
Therefore it is enough to consider only the bodies with regular semialgebraic
boundaries in $\R^n$, i.e. to assume that $\partial D$ consists of several
smooth connected components of the zero locus of a real polynomial. Also, it is
enough to prove our theorem for any connected component of $D$ separately, so
we will assume that $D$ is connected. On the other hand, for any finite $m$
there exist algebraically integrable domains in $\R^{2k}$ with $C^m$-smooth
boundaries, see \cite{Appendix}. Therefore the condition of
$C^\infty$-smoothness in Theorem \ref{mainth} cannot be reduced.

3. By an Archimedes' theorem, the volume cut off by a plane from a ball in
$\R^3$ depends polynomially on the distance of the plane from the center, i.e.
on a two-valued algebraic function in $P_3$. We need also take into account all
planes not intersecting the ball, therefore we add two single-valued functions
$P_3 \to \R^1$ equal identically to 0 and to the volume of the ball, and obtain
a four-valued algebraic function proving the algebraic integrability of the
ball in $\R^3$. Moreover, the Archimedes' theorem can be easily extended to
balls and ellipsoids in any odd-dimensional spaces. It seems likely that there
are no other examples of irreducible integrable domains with smooth boundaries
in odd-dimensional spaces, but this conjecture is not yet proved.
\end{remark}

\section{Two main examples}

\subsection{Convex domains in $\R^{2k}$}

\label{convv} Let $n$ be an even number, and $D$ a convex domain in $\R^{n}$
bounded by a compact semialgebraic variety $\partial D$ without singular
points. Choose a linear function $l: \R^{n} \to \R^1,$ whose restriction to
$\partial D$ is a Morse function (such functions exist by the Sard's lemma
applied to the Gauss map $\partial D \to \R P^{n-1}$). Let $m<M$ be both
critical values of this restriction. Denote by $A$ the complexification of
$\partial D$ (i.e. the hypersurface in $\CC^n$ distinguished by the same
polynomial equation). Also we can and will consider $l$ as the restriction to
$\R^n$ of a linear function $(\CC^n,\R^n) \to (\CC^1,\R^1)$, which we'll denote
also by $l$.

For any $t \in (m,M)$ define $V(t)$ as the volume of the set $D \cap
l^{-1}((-\infty, t])$. If $D$ is algebraically integrable, then $V$ is an
algebraic function, in particular its analytic continuation to $\CC^1$ is
finite-valued.

Further, denote by $\CC P_n$ the space of all complex affine hyperplanes in
$\CC^n$. For any $X \in \CC P_n$, consider the group
\begin{equation}
\label{cont} H_n(\CC^n, A \cup X).
\end{equation}

\begin{lemma}
There is a well-defined linear function $H_n(\CC^n, A \cup X) \to \CC^1,$ whose
value on a relative homology class is equal to the integral of the volume form
\begin{equation}
\label{vol} dx_1 \wedge \dots \wedge dx_n
\end{equation}
along an arbitrary piecewise-smooth relative cycle representing this class.
\end{lemma}

\noindent {\it Proof.} This follows from the Stokes' theorem applied to the
holomorphic form (\ref{vol}), and from the fact that the integral of this form
along any singular $n$-chain in the $(n-1)$-dimensional complex variety $A \cup
X$ is equal to zero. \hfill $\Box$ \medskip

There is a Zariski open subset $\mbox{Reg} \subset \CC P_n$ such that all
planes $X \in \mbox{Reg}$ are transversal to (maybe singular) variety $A$:
indeed, for any stratum of any algebraic Whitney stratification of $A$ the set
of planes $X$ not transversal to this stratum is a semialgebraic subvariety of
positive codimension in $\CC P_n$. Then by Thom's isotopy lemma (see e.g.
\cite{GM}) there is a locally trivial fiber bundle over $\mbox{Reg}$, whose
fiber over the point $\{X\}$ is the pair $(\CC^n, A \cup X)$. Consider the
associated homological bundle over $\mbox{Reg}$, whose fiber over $\{X\}$ is
the group (\ref{cont}).

This fiber bundle is locally trivialized, i.e. it carries the flat {\em
Gauss-Manin connection} (see e.g. \cite{AGLV}) defined by continuous shifts of
cycles into the neighboring fibers in correspondence with any local
trivialization of the initial fiber bundle of pairs $(\CC^n, A \cup X)$. This
connection is well-defined because the homological classes of these moved
cycles do not depend on the exact choice of this local trivialization. In
particular, if we fix a point $\{X_0\} \in \mbox{Reg}$ and a class $\gamma \in
H_n(\CC^n, A \cup X_0)$, then a function $\beta_{\gamma}(X)$ arises in any
simply-connected neighborhood of $\{X_0\}$ in $\mbox{Reg}$: its value at the
point $X$ is equal to the integral of the form (\ref{vol}) along the cycle
obtained from $\gamma$ by the Gauss-Manin connection over any path connecting
$X_0$ and $X$ in our neighborhood. This function is analytic in this
neighborhood, and can be continued to an analytic function on entire manifold
$\mbox{Reg}$. The ramification of this function along paths in $\mbox{Reg}$
depends on the monodromy action of the group $\pi_1(\mbox{Reg}, \{X_0\})$ on
$H_n(\CC^n, A \cup X_0)$.

\begin{figure}
\begin{center}
\begin{picture}(110,105)
\bezier{70}(60,100)(30,100)(30,80) \bezier{70}(30,80)(30,63)(50,63)
\bezier{70}(50,63)(80,63)(80,84) \bezier{70}(80,84)(80,100)(60,100)
\put(20,40){\line(1,0){90}} \put(110,40){\line(-1,-2){20}}
\put(90,0){\line(-1,0){90}} \put(0,0){\line(1,2){20}}
\put(10,20){\line(1,0){90}} \put(30,20){\circle*{2}} \put(80,20){\circle*{2}}
\put(28,14){$m$} \put(78,13){$M$} \put(39,18.5){$\blacksquare$}
\put(39.5,21){\vector(-1,0){8}} \put(31.5,19){\vector(1,0){8}}
\put(30,21){\oval(3,3)[t]} \put(30,19){\oval(3,3)[b]} \put(39,15){$\tau $}
\put(28,25){$\alpha(m)$} \put(28.5,19){\line(0,1){2}}
\put(40,64){\line(0,1){32}} \put(40.5,19){\vector(1,0){38}}
\put(78.5,21){\vector(-1,0){36.7}} \put(80,21){\oval(3,3)[t]}
\put(80,19){\oval(3,3)[b]} \put(81.5,19){\line(0,1){2}}
\put(55,25){$\alpha(M)$} \put(40,61){\vector(0,-1){28}}
\put(102,17.5){${\mathbb R}^1$} \put(100,34){${\mathbb C}^1$}
\put(31,80){\line(1,1){8.5}} \put(31.5,84){\line(1,1){8}}
\put(32.2,88){\line(1,1){7.3}} \put(31.5,76){\line(1,1){8}}
\put(32,72){\line(1,1){7.5}} \put(34,68){\line(1,1){5.2}} \put(42,52){\Large
$l$}
\end{picture}
\caption{Picard-Lefschetz monodromy of integration cycles for a convex domain
in even-dimensional space} \label{onne}
\end{center}
\end{figure}
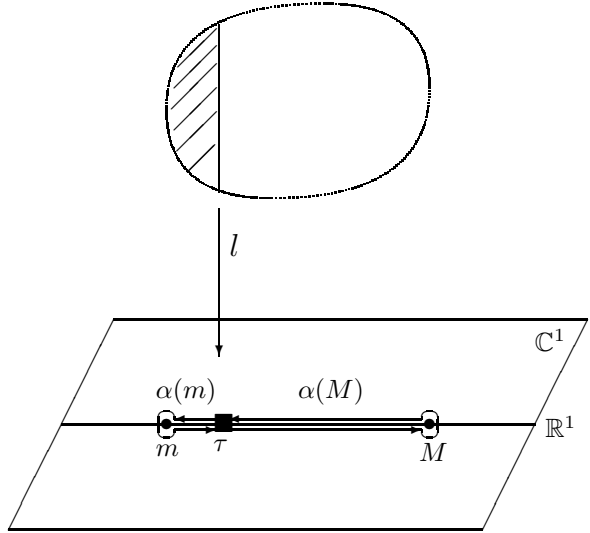

For instance, let $X_0 = X(\tau ) \in \mbox{Reg}$ be the hyperplane
$l^{-1}(\tau ) \subset \CC^n,$  $\tau  \in (m,M) \subset \R^1$, and
$\gamma=\gamma(\tau ) \in H_n(\CC^n, A \cup X(\tau ))$ the class of the figure
$D \cap l^{-1}((-\infty,\tau ])$ oriented by the form (\ref{vol}). Denote by
${\mathcal Reg}$ the set of all $t\in \CC^1$ such that $X(t) \equiv l^{-1}(t)
\in \mbox{Reg}$. For all $t \in {\mathcal Reg}$ close to $\tau $, the function
$\beta_{\gamma}( X(t))$ coincides with the volume function $V(t)$, hence their
analytic continuations to entire ${\mathcal Reg}$ also coincide. So, in order
to investigate the ramification of the volume function, let us study the orbit
of this element $\gamma$ under the monodromy action of $\pi_1({\mathcal Reg})$.

It follows easily from the Picard-Lefschetz formula (see e.g. \cite{APLT}) that
the loop $\alpha(m) \in \pi_1({\mathcal Reg})$ (see Fig. \ref{onne}) moves the
element $\gamma$ to $-\gamma$, and the loop $\alpha(M)$ moves $\gamma$ to $2
[D] - \gamma$, where $[D]$ is the homology class of the entire domain $D$.
Hence the composition of these two operators moves $\gamma$ to $\gamma + 2[D]$.
The relative cycle $[D]$ does not depend on $t$ and is invariant under the
action of the group $\pi_1({\mathcal Reg})$, hence iterating the operator
$\alpha(M) \circ \alpha(m)$ we obtain consecutively $\gamma+ 4[D]$, $\gamma +
6[D]$, etc. But the volume of $D$ is positive, hence our analytic function
takes infinitely many different values at one and the same point $\tau $, and
cannot be algebraic.

\begin{remark} \rm
Picard-Lefschetz operators act differently in spaces $\R^n$ of different
parities, because the intersection form in the middle homology group of an
$(n-2)$-dimensional complex manifold is symmetric if $n$ is even and
antisymmetric if $n$ is odd. In particular, if $n$ is odd, then the loops
$\alpha(m)$ and $\alpha(M)$ act trivially on $\gamma(t)$, which makes the
Archimedes' example possible.
\end{remark}

\subsection{General compact domain with smooth algebraic boundary in $\R^2$}

Consider first the sample domain $D$ shown in Fig. \ref{two}. Let $l:\R^2 \to
\R^1$ be the projection downwards in the page, where the target line $\R^1$ is
oriented to the right. The restriction of this function to $\partial D$ is a
strictly Morse function. Let $\tau  \in \R^1$ be a value greater than the
global minimum of $l$ on $\partial D$, but lower than all other its critical
values. For real $t \approx \tau $, let $ \gamma(t) \in H_2(\CC^2, A \cup
X(t))$ be the class of the positively oriented domain $D \cap l^{-1}((-\infty,
t])$, see Fig. \ref{two}a. Let us increase $t$ until the right-hand boundary
segment of this domain meets a critical point of $l|_{\partial D}$, shortly
before this meeting move $t$ into the complex domain in $\CC^1$ and go around
the corresponding critical value. The Gauss-Manin connection over entire this
path moves $\gamma(t)$ into the homology class in $H_2(\CC^2, A \cup X(t))$ of
the sum of two domains as shown in Fig. \ref{two}b; notice that these two
domains should be taken with opposite orientations, so that their common
boundary in $X(t)$ is a single oriented segment. Further, we decrease $t$ until
one of endpoints of this segment becomes again a critical point of
$l|_{\partial D}$, turn $t$ around the corresponding critical value (see Fig.
\ref{two}c), etc. It is important that when $t$ passes for the second time the
critical value at the local (but not global) maximum of $l$ on $\partial D$
(between Figs. \ref{two}c and \ref{two}d), the corresponding leaf of the
integral of $dx_1 \wedge dx_2$ along $\gamma(t)$ has a regular point: indeed,
the derivative of this integral over $t$ is equal to the length of the boundary
segment in $X(t)$. Therefore we do not need to make a whole circle around this
critical value: instead, we miss this critical value by a half-circle in either
half-plane of $\CC^1$ and continue to increase $t$ along the real line, see
Fig. \ref{two}d.

\unitlength 0.35mm
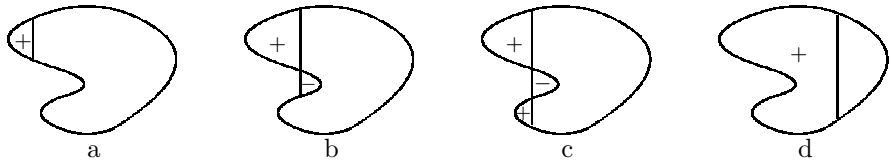
\begin{figure}
\begin{picture}(400,60)
\mbox{
\begin{picture}(80,60)
\bezier{70}(20,55)(-10,45)(20,35) \bezier{50}(20,35)(45,28)(25,23)
\bezier{50}(25,23)(10,17)(25,11) \bezier{50}(25,11)(35,7)(45,11)
\bezier{100}(45,11)(90,39)(50,55) \bezier{40}(50,55)(35,60)(20,55)
\put(14.5,37){\line(0,1){16}} \put(35,0){a} \put(7.5,42){\footnotesize $+$}
\end{picture}
} \mbox{
\begin{picture}(80,60)
\bezier{70}(20,55)(-10,45)(20,35) \bezier{50}(20,35)(45,28)(25,23)
\bezier{50}(25,23)(10,17)(25,11) \bezier{50}(25,11)(35,7)(45,11)
\bezier{100}(45,11)(90,39)(50,55) \bezier{40}(50,55)(35,60)(20,55)
\put(26,23){\line(0,1){33.3}} \put(35,0){b} \put(14,41){\footnotesize $+$}
\put(25.9,26){\scriptsize $-$}
\end{picture}
} \mbox{
\begin{picture}(80,60)
\bezier{70}(20,55)(-10,45)(20,35) \bezier{50}(20,35)(45,28)(25,23)
\bezier{50}(25,23)(10,17)(25,11) \bezier{50}(25,11)(35,7)(45,11)
\bezier{100}(45,11)(90,39)(50,55) \bezier{40}(50,55)(35,60)(20,55)
\put(24,12.5){\line(0,1){43}} \put(35,0){c} \put(14,41){\footnotesize $+$}
\put(24.8,26){\footnotesize $-$} \put(17.5,15){\scriptsize $+$}
\end{picture}
} \mbox{
\begin{picture}(80,60)
\bezier{70}(20,55)(-10,45)(20,35) \bezier{50}(20,35)(45,28)(25,23)
\bezier{50}(25,23)(10,17)(25,11) \bezier{50}(25,11)(35,7)(45,11)
\bezier{100}(45,11)(90,39)(50,55) \bezier{40}(50,55)(35,60)(20,55)
\put(50,15){\line(0,1){39}} \put(35,0){d} \put(32,37){\footnotesize $+$}
\end{picture}
}
\end{picture} \caption{Transformations of the integration cycle in $\R^2$}
\label{two}
\end{figure}

Finally our boundary segment in $X(t)$ shrinks to a point, namely to the global
maximum point of $l$ on $\partial D$. At the instant of this degeneration the
cycle obtained from $\gamma(t)$ by the Gauss-Manin connection over our path
would coincide with entire domain $D$. However we do not let $t$ to meet the
maximum value, but shortly before this rotate $t$ again in the complex domain,
and then go back along the same path. When $t$ returns to the initial position
$\tau $, the Gauss--Manin connection over the obtained loop in ${\mathcal Reg}$
moves the cycle $\gamma(\tau )$ to $2[D]-\gamma(\tau )$. The rest of the proof
repeats that in the previous example.

In the case of a general connected domain $D \subset \R^2$ with smooth
boundary, we proceed in the similar way. Namely, we take the exterior component
$\overline{\partial} D$ of $\partial D$ (which is homeomorphic to $S^1$) and a
linear function $l:(\CC^n,\R^n) \to (\CC, \R)$, whose restriction to
$\overline{\partial} D$ is Morse. Denote by $m<M$ the global extremal values of
$l$ on $\overline{\partial} D$. The corresponding global extremal points of $l$
break $\overline{\partial} D$ into two segments. Consider the set of pairs
$(a,b) \in \overline{\partial} D \times \overline{\partial} D$ such that $a$
belongs to one such part and $b$ to the other, and $l(a) = l(b)$. It is easy to
see that this set is a smooth one-dimensional submanifold with boundary, and
this boundary consists of two points $(a,a)$ where $a$ is either the global
minimum or the global maximum point of $l$ on $\partial D$. Although this
submanifold can be not connected, its boundary points should belong to one and
the same its connected component diffeomorphic to a segment. So, we start from
its point $(a,b)$ where $\{a,b\} = \partial D \cap l^{-1} (\tau)$, $\tau
\approx m$, which is close to one endpoint of this segment, and go along this
segment towards the other endpoint, watching the corresponding common values of
$l$ and making appropriate rotations in ${\mathcal Reg}$ each time when $a$ or
$b$ approaches a critical point of $l|_{\overline{\partial} D}$. This gives us
a path in ${\mathcal Reg}$. Gauss-Manin connection over this path moves the
cycle $\gamma(\tau )$ into the entire domain bounded by $\overline{\partial}
D$; in particular the value of the analytic continuation of the volume function
$V(t)$ along this path close to its endpoint is almost equal to the area of
this domain. Then we turn in $\CC^1$ around the maximal value $M$ of $l$ on
$\partial D$ and come back to $\tau $ along the same path. The monodromy along
this path in ${\mathcal Reg}$ moves the homology class $\gamma(\tau )$ into
$2[D] - \gamma(\tau )$. The rest of the construction is the same as above.

\begin{remark} \rm
To be rigorous, there can be other points of $\CC^1 \setminus {\mathcal Reg}$
on the segment $[m,M]$ apart from the critical values of the restriction of $l$
to $\overline{\partial} D$. Therefore, constructing our path in ${\mathcal
Reg}$, we need to avoid these points along small arcs in the complex domain.
This does not affect our consideration, because our integration cycle (and its
area function) behave regularly when $t$ passes these additional points.
\end{remark}

\section{Proof of the main theorem}

\subsection{On finite reflection groups}
We use the following well-known facts about finite reflection groups, see e.g.
\cite{bour}, \cite{AGLV}.

\begin{prop}
\label{refle} Let $\Z^N$ be an integer lattice with integer-valued symmetric
bilinear form $\langle \cdot , \cdot \rangle$; let $\{e_j\} \subset \Z^N$ be a
finite collection of elements of length $\sqrt{2}$ $($i.e. $\langle e_j, e_j
\rangle = 2$ for each $j)$ generating the entire $\Z^N$ as a $\Z$-module. Let
$G$ be the subgroup in $\mbox{SL}(N, \Z)$ generated by reflections
corresponding to these elements $e_j$, acting on $\Z^N$ by the formula
\begin{equation}
\label{reflect} R_j : a \mapsto a-\langle e_j, a \rangle e_j.
\end{equation}
Suppose that the orbits in $\Z^N$ of all generating elements $e_j$ under the
action of the group $G$ are finite. Then

1. The group $G$ is finite.

2. The form $\langle \cdot, \cdot \rangle$ is non-degenerate\footnote{and even
positive definite, which is less important for us now} : if $\langle e_j , a
\rangle = 0$ for any $j$, then $a=0$.
\end{prop}

\subsection{Reflection group related with a smooth semialgebraic domain in $\R^{2k}$}

Suppose that $n$ is even. Let $D$ be a bounded connected domain with
semialgebraic non-singular boundary in $\R^n$. This boundary can consist of
several connected components, but exactly one of them separates the entire $D$
and all other components of $\partial D$ from the infinity. Let us denote this
component by $\overline{\partial} D$. Let $A$ be the complexification of
$\overline{\partial} D$: this is a hypersurface in $\CC^n$, which can have
singularities apart from a neighborhood of $D$. Let us fix some complex
semialgebraic Whitney stratification of $A$ (see e.g. \cite{GM}), all whose
strata of dimension $<n-1$ do not meet $\partial D$.

\begin{definition} \rm
\label{gener} A linear function $l:(\CC^n,\R^n) \to (\CC^1,\R^1)$ with real
coefficients is {\em generic} with respect to the hypersurface $A$ if

1. Its restriction to $\partial D$ is a strictly Morse function;

2. For any critical value $t \in \R^1$ of this restriction, the corresponding
complex hyperplane $X(t) = l^{-1}(t) \subset \CC^n$ is transversal to all
strata of $A$ apart from the corresponding critical point of this restriction.
\end{definition}

Generic linear functions are dense in the space of all linear functions; let
$l$ be one of them. Recall the notation ${\mathcal Reg}$ for the set of all
values $t \in \R^1$ such that the plane $X(t)$ is transversal to $A$.

Let $a_j \in \overline{\partial} D$ be a critical point of
$l|_{\overline{\partial} D}$, and $t_j$ the corresponding critical value
$l(a_j)$. For any sufficiently small $\varepsilon>0$ and any $\nu>0$
sufficiently small with respect to $\varepsilon$, denote by $B_\varepsilon
(a_j)$ the ball of radius $\varepsilon$ in $\CC^n$ centered at $a_j$, and
suppose that $t \in \CC^1$ is an arbitrary point in the punctured
$\nu$-neighborhood of $t_j$. Then by Milnor's theorem (see e.g. \cite{AGLV})
the (reduced modulo a point) homology group $\tilde H_{n-2}(B_\varepsilon (a_j)
\cap A \cap X(t))$ is isomorphic to $\Z$. The chain of isomorphisms
\begin{equation}
\label{redu1}
\begin{array}{l} H_n(B_\varepsilon (a_j), B_\varepsilon(a_j) \cap (A
\cup X(t)))
\to H_{n-1}(B_\varepsilon (a_j) \cap (A \cup X(t))) \to \\
\to H_{n-1}(B_\varepsilon (a_j) \cap X(t), B_\varepsilon (a_j) \cap A \cap
X(t)) \to \tilde H_{n-2}(B_\varepsilon (a_j) \cap A \cap X(t))
\end{array}
\end{equation} (in which the first and the third arrows are boundary
homomorphisms, and the second one is a fragment of the exact sequence of the
pair $(B_\varepsilon(a_j) \cap (A \cup X(t)), B_\varepsilon(a_j) \cap A)$)
proves the same for its left-hand group
\begin{equation}
\label{local} H_n(B_\varepsilon (a_j), B_\varepsilon (a_j) \cap (A \cup X(t))).
\end{equation}

The monodromy action on these groups, defined by the rotation of $t$ around the
value $t_j$ in its $\nu$-neighborhood, commutes with all these isomorphisms and
takes one generator of this group (\ref{local}) into the other (i.e. it is
multiplication by $-1$).

Let us choose an arbitrary such generator $\delta_j$ for $t = t_j+ \nu/2$,
extend this choice by the Gauss-Manin connection to similar groups
(\ref{local}) for all $t \in (t_j, t_j+\nu)$, and define the function
$\varphi_j: (t_j, t_j+\nu) \to \CC^1$, whose value $\varphi_j(t)$ is equal to
the integral of the form (\ref{vol}) along this generator of (\ref{local}). It
is easy to calculate that this function is real or purely imaginary depending
on the parity of the Morse index of $l|_{\partial D}$ at $a_j$, and its
absolute value vanishes asymptotically as $(t-t_j)^{(n+1)/2}$ when $t$ tends to
$t_j$. Also, it is analytic and can be extended to a neighborhood of the point
$t_j + \nu/2$ in $\CC^1$. Let again $[m,M]$ be the segment of values
$l(\overline{\partial} D) \subset \R^1.$ Since the set $\CC^1 \setminus
{\mathcal Reg}$ is finite, there is a neighborhood $U([m,M])$ of this segment,
such that all points of this set in $U([m,M])$ are real; they include the set
$\Sigma$ of critical values of $l|_{\overline{\partial} D}$. We can and will
assume that $0 \in [m,M] \cap {\mathcal Reg}$. Let us connect the point $0$
with all points $t_j+ \nu/2$ by arbitrary paths in the domain $U([m,M]) \cap
{\mathcal Reg} \cap \{t: Im(t) > 0\}$: \unitlength 1.5mm \linethickness{0.6pt}
\begin{picture}(40,3)
\put(0,0){\line(1,0){40}} \put(19.3,-0.8){\Large $0$} \put(37,1){$\R$}
\put(2,0){\circle*{1.2}} \put(12,0){\circle*{1.2}} \put(33,0){\circle*{1.2}}
\put(11.75,0){\oval(16.5,5)[t]} \put(16.75,0){\oval(6.5,3)[t]}
\put(27.25,0){\oval(14.5,5)[t]} \put(2,0){\oval(3,3)[t]}
\put(12,0){\oval(3,3)[t]} \put(33,0){\oval(3,3)[t]}
\end{picture} \label{pic3} . Define the function germs
$\Phi_j: (\CC^1,0) \to \CC$ at the point $0$ as analytic continuations of
functions $\varphi_j$ along these paths. Their values at $0$ are equal to the
integrals of (\ref{vol}) along the relative cycles $\tilde \Delta_j$ in $\CC^n
\mbox{ mod } (A \cup X(0))$ obtained from the chosen generators of
(\ref{local}) by the Gauss-Manin connection over these paths. If we realise
this Gauss-Manin connection by the local trivialization of the fiber bundle of
pairs $(\CC^n, A \cup X(t))$ following from the Thom's isotopy lemma, then the
closures of these cycles $\tilde \Delta_j$ do not meet the singular locus of
$A$. Denote by ${\stackrel{\circ}{A}}$ the non-singular part of $A$, then the
cycles obtained in this way define certain elements $ \Delta_j \in H_n(\CC^n
\setminus \mbox{sing}(A), {\stackrel{\circ}{A}} \cup X(0))$.

Let $\mathfrak{F}$ be the subgroup of the space of all complex-valued function
germs at $0$, consisting of integer linear combinations of our function germs
$\Phi_j$. This subgroup is the image of a homomorphism of an integer lattice
(consisting of all integer linear combinations of symbols $\Phi_j$) into the
space of germs, hence also an integer lattice. For any $i, j$ define the scalar
product $\langle \Phi_i, \Phi_j \rangle$ as follows: we take two
$(n-2)$-dimensional cycles in ${\stackrel{\circ}{A}} \cap X(0)$ obtained from
$\Delta_i$ and $\Delta_j$ by the composite homomorphism
\begin{equation}
\label{reduc}
\begin{array}{r}
H_n(\CC^n \setminus \mbox{sing}(A), {\stackrel{\circ}{A}} \cup X(0)) \to
H_{n-1}({\stackrel{\circ}{A}} \cup X(0)) \to \\ \to H_{n-1}(X(0),
{\stackrel{\circ}{A}} \cap X(0)) \to H_{n-2}({\stackrel{\circ}{A}} \cap X(0))
\end{array}
\end{equation}
similar to (\ref{redu1}), calculate the intersection index of these
$(n-2)$-dimensional cycles in the complex variety ${\stackrel{\circ}{A}} \cap
X(0),$ and multiply this intersection index by $(-1)^{1+n/2}.$ Since $n$ is
even, $\langle \Phi_i, \Phi_j \rangle = \langle \Phi_j, \Phi_i \rangle$ for any
$i, j$. The scalar squares $\langle \Phi_j, \Phi_j \rangle$ of all generating
elements $\Phi_j$ are equal to $2$ (see e.g. \cite{AGLV}, \S II.1.3).

\begin{lemma}
\label{lem1} The scalar products $\langle \Phi_i, \Phi_j \rangle$ can be
extended by linearity $($in a unique way$)$ to a symmetric bilinear form on the
lattice $\mathfrak{F}$.
\end{lemma}

\noindent {\it Proof.} Suppose that the linear combination $c_1 \Phi_1 + \dots
+ c_r \Phi_r$ defines the identically zero function germ. Consider the loop in
$U([m,M]) \cap {\mathcal Reg}$, consisting of the distinguished path from $0$
to $t_j +\nu/2$, rotation around $t_j$ along the circle of radius $\nu/2$, and
return to $0$ by our distinguished path. By the Picard-Lefschetz formula, the
{\it analytic continuation of our zero germ} along this path turns it into the
germ $\Phi_j$ taken with the coefficient $\pm \langle (c_1 \Phi_1 + \dots + c_r
\Phi_r), \Phi_j \rangle.$ The germ $\Phi_j$ is not equal to identical zero,
hence this coefficient should vanish for any $j$. \hfill $\Box$ \medskip

Each basic element $\Phi_j$ defines a reflection in the space $\mathfrak{F}$
acting by
\begin{equation}
\label{PL2} \Phi \mapsto \Phi - \langle \Phi_j, \Phi \rangle \Phi_j;
\end{equation}
by Picard--Lefschetz formula it describes the analytic continuation of an
arbitrary germ $\Phi \in \mathfrak{F}$ along the loop constructed in the proof
of lemma \ref{lem1}. These reflections preserve the scalar product $\langle
\cdot , \cdot \rangle$ and thus generate a subgroup of the group of
automorphisms of $(\mathfrak{F}; \langle \cdot , \cdot \rangle$).

\begin{prop}
\label{pro3} If the domain $D$ is algebraically integrable, then the subgroup
generated by operators $($\ref{PL2}$)$ in the group of automorphisms of the
lattice $\mathfrak{F}$ is finite.
\end{prop}

\noindent {\it Proof}. The set of all complex hyperplanes $X(t) = l^{-1}(t),$
$t \in \CC^1$, forms a line in the space $\CC P_n$. We can assume that $\nu$ is
small enough, so that $(m,m+\nu) \subset {\mathcal Reg}$. Let $\tau =m+\nu/2$,
and denote by $\gamma(X(\tau ))$ the class in $H_n(\CC^n, A \cup X(\tau ))$ of
the positively oriented domain $D \cap l^{-1}((-\infty, \tau ])$. Consider also
the analytic function $\mbox{Vol}$ on $\mbox{Reg} \subset \CC P_n$ defined by
integrals of the form (\ref{vol}) along the similar cycles $\gamma(X) \in
H_n(\CC^n, A \cup X)$ obtained from $\gamma(X(\tau ))$ by Gauss-Manin
connection over different paths in $\mbox{Reg}$ connecting $X(\tau )$ and $X$.

\begin{lemma}
\label{lem2} The set of restrictions of different leaves of the function
$\mbox{\rm Vol}$ to neighborhoods of the point $0$ in the line $\CC^1 \subset
\CC P_n$ of planes $X(t)$ contains all function germs $\Phi_j: (\CC^1, 0) \to
\CC^1$ corresponding to all critical points of the restriction of $l$ to
$\overline{\partial} D$.
\end{lemma}

\noindent {\it Proof of lemma \ref{lem2}}. It is enough to prove that for any
$j$ the integral functions $\varphi_j:(t_j, t_j+\nu) \to \CC^1$ defined above
appears among restrictions of the function $\mbox{Vol}$ to the interval
$(t_j,t_j+\nu) \subset \CC^1$. Consider the set of all non-singular points of
the hypersurface $A$, at which the second fundamental form of this hypersurface
is non-degenerate. Its complement has complex codimension 1 in $A$, hence this
set is path-connected within any irreducible component of $A$. Since
$l|_{\overline{\partial} D}$ is a Morse function, this set contains all its
critical points (and all of them obviously belong to one and the same component
of $A$). Choose a smooth path in this set, connecting the critical points with
values $m$ and $t_j$. The tangent hyperplanes at these points define a path in
$\CC P_n$. Take a path in the set $\mbox{Reg} \subset \CC P_n$ escorting this
one in its thin neighborhood and connecting the planes $X(\tau )$ and
$X(t_j+\nu/2)$ in such a way that for any point $\{X\}$ of this path the
element of the group (\ref{cont}) obtained from $\gamma$ by the Gauss-Manin
connection over this path is realised by a cycle inside a small ball centered
at the neighboring tangency point of $A$ and a hyperplane parallel to $X$, in
particular this cycle for the endpoint of this path generates the group
(\ref{local}). This terminal cycle can coincide with the integration cycle
defining the function $\varphi_j$ or be opposite to it. In the first case
$\varphi_j$ is defined by the analytic continuation of the function
$\mbox{Vol}$ along this path, in the second one we need additionally go once
around the critical value $t_j$ in the set ${\mathcal Reg}$ of generic
hyperplanes $X(t)$. \hfill $\Box$

\begin{remark} \rm
In this lemma, we do not state that any germ $\Phi_j$ can be obtained from the
initial function at the point $\tau $ by the analytic continuation along a path
inside our line $\CC^1$. Moreover, the reflection group on $\mathfrak{F}$
generated by operators (\ref{PL2}) (defined by paths inside this domain only)
can be reducible, see subsection \ref{redred} below.
\end{remark}

\begin{remark} \rm
The path in $\CC P_n$ used in this proof can be chosen in an arbitrarily small
neighborhood of the set of $($complexifications of$)$ real hyperplanes. Indeed,
let us connect our critical points of $l|_{\overline{\partial} D}$ by a generic
path inside our component of $\overline{\partial} D$; then to avoid the set of
parabolic points of $A$ we can move this path slightly into the complex domain
in arbitrarily small neighbourhoods of points at which it crosses this set.
\end{remark}

Further, the set of analytic continuations of all function germs $\Phi_j$
contains all germs obtained from them by reflections (\ref{PL2}). Therefore if
the domain $D$ is algebraically integrable then the orbit of any element
$\Phi_j$ under the reflection group generated by operators (\ref{PL2}) is
finite. By proposition \ref{refle} the entire this reflection group is then
finite, in particular the symmetric form $\langle \cdot, \cdot \rangle$ on
$\mathfrak{F}$ is non-degenerate. Proposition \ref{pro3} is proved. \hfill
$\Box$ \medskip

\begin{prop}
\label{prp4} The class $[\overline{D}] \in H_n(\CC^n, A \cup X(0))$ of the
entire domain in $\R^n$ bounded by $\overline{\partial} D$ is equal to the sum
of classes $\Delta_j$ $($taken with appropriate signs$)$ over all critical
points of $l|_{\overline{\partial} D}$. In particular, the corresponding sum of
function germs $\pm \Phi_j$ is the constant function equal identically to the
volume bounded by $\overline{\partial} D$.
\end{prop}

This fact follows immediately from the following one (Lemma \ref{next} below).
For any critical value $t_j$ of $l|_{\overline{\partial} D}$ we have two
elements in the group $H_n(\CC^n, A \cup X(t_j+\nu/2))$: one (let us call it
$\nabla_+$) given by the positively oriented part of the half-space $\R^n\cap
l^{-1}((-\infty, t_j+ \nu/2])$ bounded by $\overline{\partial} D$, and the
other, $J(\nabla_-),$ obtained from the similar class $\nabla_- \in H_n(\CC^n,
A \cup X(t_j-\nu/2))$ by the Gauss--Manin connection over the arc of radius
$\nu/2$ in the upper half-plane of $\CC^1$.

Consider also absolute homology classes $\Pi_{\pm} \in H_{n-2}(A \cap X(t_j\pm
\nu/2))$ represented by naturally oriented real manifolds $\overline{\partial}
D \cap X(t_j\pm \nu/2)$ respectively. These classes are obtained from
$\nabla_\pm$ by maps similar to the composition (\ref{reduc}).

\begin{lemma}[see \cite{APLT}, Lemma 3.3 on page 121]
\label{next} 1. The difference of two relative homology classes $\nabla_+$ and
$J(\nabla_-)$ in $H_n(\CC^n, A \cup X(t_j+\nu/2))$ is equal to the image of the
vanishing cycle $\delta_j \in H_n(B_\varepsilon(a_j), A \cup X(t_j+\nu/2))$
under the identical embedding $B_\varepsilon(a_j) \hookrightarrow \CC^n$.

2. The difference in the group $H_{n-2}(A \cap X(t_j+\nu/2))$ of the class
$\Pi_+$ and the element obtained from $\Pi_-$ by the Gauss--Manin connection
over the arc of radius $\nu/2$ in the upper half-plane of $\CC^1$ can be
realized by a vanishing cycle generating the group $H_{n-2}(B_\varepsilon(a_j)
\cap A \cap X(t_j+ \nu/2)) \simeq \Z.$ \hfill $\Box$ \ $\Box$
\end{lemma}

\begin{remark} \rm
The earliest (very technical) proof of item 2 known to me is given in
\cite{Leray}; for another proof, deducing it directly from the Picard-Lefschetz
formula, see \cite{Petr} and \S V.3 of \cite{APLT}.
\end{remark}

Since the volume bounded by $\overline{\partial} D$ is positive, the sum
$\sum_j \pm \Phi_j$ mentioned in proposition \ref{prp4} is a non-zero element
of the lattice $\mathfrak{F}$. On the other hand, the image of the
corresponding cycle $[\overline{D}]$ under the map (\ref{reduc}) is equal to
zero, because the image of this cycle under the first arrow in (\ref{reduc}) is
a cycle in $A$, and the second arrow in (\ref{reduc}) is the reduction modulo
$A$. Therefore this non-zero element of $\mathfrak{F}$ belongs to the kernel of
the scalar product $\langle \cdot, \cdot \rangle$. Thus by proposition
\ref{refle} the group $G$ cannot be finite. Theorem \ref{mainth} is completely
proved. \hfill $\Box$

\subsection{Reducibility of the reflection group}
\label{redred}

The reflection group
\begin{equation}
\label{genref} (\mathfrak{F}; \langle \cdot , \cdot \rangle; \{R_j\})
\end{equation}
can be reducible, for instance this is the case if $D$ is a thin tubular
neighborhood of the standard circle embedded into $\R^4$. In this case we can
choose the linear function $l: \R^4 \to \R$ in such a way that its restriction
to $\partial D$ has four critical points $a_j$, $j = 1, 2, 3, 4$ with Morse
indices $0, 1, 2$ and $3$ respectively, and the critical values $t_2$ and $t_3$
at the points $a_2$ and $a_3$ almost coincide. Then we have the following
calculation.

\begin{prop}
For appropriate choice of orientations of vanishing cycles $\beta_i \in
H_2(B_\varepsilon (a_j) \cap A \cap X(t_j+\nu/2)) \simeq \Z,$ the intersection
matrix in $H_{2}(A \cap X(0))$ of cycles obtained from them by the Gauss-Manin
connection over any paths in the upper half of $U([m,M]) \cap {\mathcal Reg}$
$($see picture in page \pageref{pic3}$)$ is equal to
$$\begin{vmatrix}
\ 2 & -2 & \ 0 & \ 0 \\ -2 & \ 2 & \ 0 & \ 0 \\ \ 0 & \ 0 & \ 2 & -2 \\ \ 0 & \
0 & -2 & \ 2
\end{vmatrix}\ .
$$

In particular, in this case the entire action of the reflection group
$($\ref{genref}$)$ splits into the direct sum of two-dimensional ones, each of
which is isomorphic to one arising in the example of subsection \ref{convv}.
\end{prop}

\noindent {\it Sketch of the proof}. $\langle \beta_2, \beta_3 \rangle = 0$
because the corresponding critical points are distant, but their critical
values can be made arbitrarily close to one another by bending the function $l$
(which does not affect the intersection index). Further, we can assume that $0$
is slightly above the critical values $t_2$ and $t_3$. In this case the cycle
in the variety $A \cap X(0)$ obtained by the Gauss--Manin connection from
$\beta_4$ is presented by the set $\partial D \cap X(0)$ of all its real
points. Its intersection with the vanishing cycles $\beta_3$ and $\beta_2$ can
be calculated in local terms of Morse critical points $a_2$ and $a_3$ (see
\cite{GZ}, \cite{AC}, \cite{APLT}). Replacing $0$ by a value slightly below
$t_2$ and $t_3$ (which does not affect the intersection indices of cycles
obtained by Gauss--Manin connection in the upper half-plane) we obtain similar
indices for cycles obtained from $\beta_1, \beta_2$ and $\beta_3$.

Finally, we know from lemma \ref{next} that the cycle of real points in $A \cap
X(t_4 - \nu/2)$ can be realised as the sum of three cycles obtained by
Gauss-Manin connection in the upper half-plane from $\beta_1$, $\beta_2$ and
$\beta_3$. Hence, by the previous calculations, the cycle obtained from
$\beta_4$ has with this sum the same intersection index, as with $\beta_3$
only. This gives us zero for the remaining corner elements $(1,4)=(4,1)$ of the
intersection matrix. \hfill $\Box$ \medskip

In general, if the reflection group (\ref{genref}) is reducible, then the set
of all cycles $\Delta_j$ splits into collections of cycles generating germs
$\Phi_j$ that belong to different irreducible components. It is easy to see
that the sum of all function germs corresponding to all cycles $\Delta_j$ from
any such collection belongs to the kernel of the form $\langle \cdot, \cdot
\rangle$. Indeed, this sum for one such collection is obviously orthogonal to
all generators $\Phi_j$ from the other collections. On the other hand, this sum
is equal to the difference of the sum of all $\Phi_j$ over all collections and
the sum of such sums over all collections except for this one; both sums are
orthogonal to all $\Phi_j$ from our collection.

\end{document}